# On predictive probability matching priors

**Trevor J. Sweeting**[1]

*University College London*

**Abstract:** We revisit the question of priors that achieve approximate matching of Bayesian and frequentist predictive probabilities. Such priors may be thought of as providing frequentist calibration of Bayesian prediction or simply as devices for producing frequentist prediction regions. Here we analyse the $O(n^{-1})$ term in the expansion of the coverage probability of a Bayesian prediction region, as derived in [*Ann. Statist.* **28** (2000) 1414–1426]. Unlike the situation for parametric matching, asymptotic predictive matching priors may depend on the level $\alpha$. We investigate *uniformly predictive matching priors* (UPMPs); that is, priors for which this $O(n^{-1})$ term is zero for all $\alpha$. It was shown in [*Ann. Statist.* **28** (2000) 1414–1426] that, in the case of quantile matching and a scalar parameter, if such a prior exists then it must be Jeffreys' prior. In the present article we investigate UPMPs in the multiparameter case and present some general results about the form, and uniqueness or otherwise, of UPMPs for both quantile and highest predictive density matching.

## Contents

## 1. Introduction

Prior distributions that match posterior predictive probabilities with the corresponding frequentist probabilities are attractive when a major goal of a statistical analysis is the construction of prediction regions. Such priors provide calibration of Bayesian prediction or may be viewed as a Bayesian mechanism for producing frequentist prediction intervals.

It is known that exact predictive probability matching is possible in cases in which there exists a suitable transformation group associated with the model. The general group structure for parametric models starts with a group of transformations on the sample space under which the statistical problem is invariant. This group of transformations then gives rise to a group $G$ of transformations on the parameter space. From an "objective Bayes" point of view, it makes sense to choose a prior distribution that is (relatively) invariant under this group. In particular, this will ensure that initial transformation of the data will make no difference to predictive inferences. The two fundamental invariant measures on the group $G$ are the left and right Haar measures. The left (right) Haar measure is the unique left-(right-)translation invariant measure on $G$, up to a positive multiplicative constant. These measures give rise to invariant left and right Haar priors on the parameter space. In the decision-theoretic development, under suitable conditions it turns out that the right Haar prior gives rise to optimal invariant decision rules for invariant decision problems; see, for example, [1]. The left Haar prior, however, which coincides with Jeffreys' invariant prior, often gives inadmissible rules in multiparameter cases. These facts provide strong motivation for the use of the right Haar prior. In relation to predictive inference, following earlier work in [10] and [11] this intuition was further reinforced in [13], where it was shown that if such a group structure exists then the associated right Haar prior gives rise to exact predictive matching for all invariant prediction regions. Thus the predictive matching problem is solved for models that possess a suitable group structure when the prediction region is invariant.

[1]Department of Statistical Science, University College London, Gower Street, London WC1E 6BT, UK, e-mail: trevor@stats.ucl.ac.uk

When exact matching is not possible one can instead resort to asymptotic approximation and investigate approximate predictive matching. This question was explored in [4] for the case of $n$ independent and identically distributed (i.i.d.) observations. For regular parametric families the difference between the frequentist and posterior predictive probabilities is $O(n^{-1})$ and a concise expression for this difference was obtained in [4] by using the auxiliary prior device introduced by P. J. Bickel and J. K. Ghosh in [3]. This technical device has proved to be extremely valuable for the theoretical comparison of Bayesian and frequentist inference statements, or simply as a Bayesian device for obtaining frequentist results. It has been particularly useful for deriving probability matching priors (see, for example, [8], [9] and the review in [5]) and for studying properties of sequential tests ([16]).

In order to find an approximate predictive probability matching prior, one sets the $O(n^{-1})$ discrepancy to zero and attempts to solve the resulting partial differential equation (PDE); a number of examples are given in [4]. We briefly review the main results in [4] in Section 2. Two main issues arise from this analysis. Firstly, the PDE for a predictive matching prior may be difficult to solve analytically. The second, and more fundamental, issue is that, except in special cases, the resulting matching prior will depend on the desired predictive probability level $\alpha$. If there does exist a prior that gives rise to predictive probability matching for all $\alpha$ then we shall refer to it as a *uniformly predictive matching prior* (UPMP). Of course, in the case of a transformation model and an invariant prediction region we already know from [13] that the right Haar prior must be a solution of the PDE. It is instructive to demonstrate this directly and this is done in the Appendix for quantile matching. Since the definition of the right Haar prior depends on a specific group of transformations on the parameter space, we need to study the effect of parameter transformation on the quantities appearing in the PDE. For this reason, it is natural to regard Fisher's information, $g$, as a Riemannian metric tensor so that transformational properties of $g$ and the other quantities that appear in the PDE can be studied.

In the case of quantile matching and a single real parameter, it has already been shown in [4] that if there exists a UPMP then it must be Jeffreys' invariant prior. This result therefore extends the exact Haar prior result for transformation models to the most general models for which approximate uniform matching is possible. However, it is clear from examples discussed in [4] and from the general theory for transformation models in [13] that this result will not hold in the multiparameter case. For example, the unique UPMP for the normal model with unknown mean and variance is the right Haar prior, or Jeffreys' independence prior, whereas Jeffreys' prior is the left Haar prior.

The main purpose of the present article is to investigate the general form of UPMPs whenever they exist. In particular, we explore the uniqueness or otherwise of the right Haar prior as a UPMP for quantile matching in the case of a transformation model. Although UPMPs can exist outside of transformation models, such situations would seem to occur rarely. The main results are given in Section 3 and 4. In Section 3 we explore the form of the UPMP for quantile matching. In addition to confirming that the right Haar prior is a UPMP in suitable transformation models, as discussed above, we obtain the general form of the UPMP whenever one exists and show that this prior is unique. In particular, it follows that for transformation models there are no priors other than the right Haar prior that give approximate uniform predictive quantile matching. In Section 4 we consider probability matching based on highest predictive density regions, which are particularly relevant for multivariate data. The scalar parameter case is clear-cut and was essentially treated in [4], where it was shown that if there exists a UPMP then it is unique. However, unlike quantile matching, this UPMP is not necessarily Jeffreys' prior. The situation is less straightforward in the multiparameter case. We show that, under a certain condition, if there exists a UPMP then it is unique. If this condition is not satisfied then either there will be no UPMP or there will exist an infinite number of UPMPs. This section provides predictive versions of results for highest posterior density regions obtained by J. K. Ghosh and R. Mukerjee in [8] and [9]. We end with some discussion in Section 5.

## 2. Review of predictive probability matching priors

We begin by introducing the notation and reviewing the main results in [4] on predictive probability matching priors. We consider only the case of i.i.d. observations in this article, but the results would be expected to hold more generally under suitable conditions. Suppose then that $X_1, X_2, \ldots$ is a sequence of independent observations having the same distribution as the (possibly vector-valued) continuous random variable $X$ with density $f(\cdot; \theta)$, where $\theta = (\theta_1, \ldots, \theta_p) \in \Omega$ is an unknown parameter and $\Omega$ is an open subset of $\Re^p$. Consider the problem of predicting the next observation, $X_{n+1}$, based on the first $n$ observations, $d = (X_1, X_2, \ldots, X_n)$. We assume regularity conditions on $f$ and $\pi$, as detailed in [4]. In particular, the support of $X$ is assumed to be independent of $\theta$.

Consider first the case of univariate $X$. Let $q(\pi, \alpha, d)$ denote the $1 - \alpha$ quantile of the posterior predictive distribution of $X_{n+1}$ under the prior $\pi$. That is, $q(\pi, \alpha, d)$ satisfies the equation

$$(2.1) \qquad P^\pi(X_{n+1} > q(\pi, \alpha, d) | d) = \alpha.$$

Let $q(\theta, \alpha)$ be the $1 - \alpha$ quantile of $f(\cdot; \theta)$; that is

$$(2.2) \qquad \int_{q(\theta,\alpha)}^{\infty} f(u; \theta) du = \alpha.$$

Write $\partial_t = \partial/\partial \theta_t$ and let $f_t(u; \theta) = \partial_t f(u; \theta)$. Define

$$(2.3) \qquad \mu_t(\theta, \alpha) = \int_{q(\theta,\alpha)}^{\infty} f_t(u; \theta) du.$$

Finally, let $g(\theta)$ be the per observation Fisher information matrix, which we assume to be non-singular for all $\theta \in \Omega$, and let $g_{st}$ and $g^{st}$ be the $(s,t)$th elements of $g$ and $g^{-1}$ respectively.

Using the approach of [3] and [7] in which an auxiliary prior is introduced and finally allowed to converge weakly to the degenerate measure at $\theta$, it follows from equations (3.3) and (3.4) in [4] that

$$(2.4) \qquad P_\theta(X_{n+1} > q(\pi, \alpha, d)) = \alpha - \frac{\partial_s\{g^{st}(\theta)\mu_t(\theta,\alpha)\pi(\theta)\}}{n\pi(\theta)} + o(n^{-1}).$$

Here and elsewhere we use the summation convention. We will say that $\pi$ is a *level-$\alpha$ predictive probability matching prior* if it satisfies the equation

$$(2.5) \qquad \partial_s\{g^{st}(\theta)\mu_t(\theta,\alpha)\pi(\theta)\} = 0.$$

From (2.4), such a prior $\pi$ matches the Bayesian and frequentist predictive probabilities to $o(n^{-1})$. Clearly, in general a solution of (2.5) will depend on the particular level $\alpha$ chosen. This is demonstrated in [4] for the specific example in which the observations are from a $N(\theta, \theta)$ distribution. Recalling the discussion in Section 1, we refer to a prior for which (2.5) holds for all $\alpha$ as a *uniformly predictive matching prior* (UPMP). In the case $p=1$, it was shown in [4] that if there exists a UPMP then this prior must be Jeffreys' prior. As noted in [4], when no UPMP exists then the formula on the left-hand side of (2.5) may still be useful for comparing alternative priors.

Moving to the multiparameter case, examples in [4] illustrate that the above result on Jeffreys' prior no longer holds. An illustration of this is Example 2 in [4], which is the location-scale model $f(x;\theta) = \theta_2^{-1} f^*(\theta_2^{-1}(x-\theta_1))$. In this case there exists a UPMP given by $\pi(\theta) \propto \theta_2^{-1}$, which is the right Haar prior for this model under the location-scale transformation group, whereas Jeffreys' prior is the left-invariant prior $\pi(\theta) \propto \theta_2^{-2}$ under this group.

In the case where $X$ is possibly vector-valued, the coverage properties of highest predictive density regions are investigated in [4]. This investigation mirrors that in [8] and [9] for highest posterior density regions. Let $m(\theta, \alpha)$ be such that

$$\int_A f(u;\theta)du = \alpha,$$

where $A = A(\theta, \alpha) = \{u : f(u;\theta) \geq m(\theta, \alpha)\}$ and define

$$\xi_j(\theta, \alpha) = \int_A f_j(u;\theta)du.$$

Let $H(\pi, \alpha, d)$ be the level-$\alpha$ highest predictive density region under the prior $\pi$. Then, as for quantile matching, it follows from the results in Section 5 of [4] that

$$P_\theta(X_{n+1} \in H(\pi, \alpha, d)) = \alpha - \frac{\partial_s\{g^{st}(\theta)\xi_t(\theta,\alpha)\pi(\theta)\}}{n\pi(\theta)} + o(n^{-1}).$$

Thus $\pi$ is a level-$\alpha$ predictive probability matching prior if and only if it satisfies the equation

$$(2.6) \qquad \partial_s\{g^{st}(\theta)\xi_t(\theta,\alpha)\pi(\theta)\} = 0.$$

Once again we see that in general the solution $\pi$ will depend on the level $\alpha$. Examples are given in [4] in which there are no priors that satisfy (2.6) for all $\alpha$. Moreover, even in the case $p=1$, if there does exist a unique prior satisfying (2.6) for all $\alpha$ then it is not necessarily Jeffreys' prior.

## 3. UPMPs: quantile matching

As discussed in Section 2 we know that when $p = 1$ and a UPMP exists for quantile matching as in (2.1), then it must be Jeffreys' prior. However, it need not be Jeffreys' prior when $p > 1$. Under a suitable group structure on the model, the results in [13] imply that the associated right Haar prior gives exact predictive matching, since the prediction region here is invariant. Thus in these cases the right Haar prior must also be a solution of equation (2.5). It is instructive to demonstrate directly that this is indeed the case.

First note from the product rule that equation (2.5) is equivalent to

$$(3.1) \qquad g^{st}(\theta)\mu_t(\theta,\alpha)\partial_s\lambda(\theta) + \partial_s\{g^{st}(\theta)\mu_t(\theta,\alpha)\} = 0,$$

where $\lambda(\theta) = \log \pi(\theta)$. Suppose that there exists a group $G$ of bijective transformations on the sample space under which the statistical problem is invariant. Further assume, as in [13], that $G = \Omega$, a locally compact topological group. In this case the distribution of $X$ under $\theta$ is the same as that of $\theta X$ under $e$, the identity element of the group, with $\theta$ regarded as an element of the transformation group $G$. Then there exist unique (up to a multiplicative constant) left-invariant and right-invariant Haar measures on $G$, giving left and right Haar priors on the parameter space. In the following we denote the right Haar prior density on $\Omega$ by $\pi^H$. The proof of the following theorem is given in the Appendix.

**Theorem 3.1.** *Under the above group structure the right Haar prior satisfies equation (3.1).*

Two questions naturally arise. First, if the above group structure exists then can there be UPMPs other than the right Haar prior? The answer to this question turns out to be "no," as follows from Theorem 3.2 below. Second, if the above group structure does not exist can there still be a UPMP? The answer to this question is "yes." An example in the case $p = 1$ is given in Section 3 of [4] for which there is no suitable group structure but there is still a unique UPMP, which must of course be Jeffreys' prior.

We now establish the general form of the UPMP whenever it exists and show that it is unique. This is a multiparameter version of Theorem 1 in [4]. Let $F(x;\theta)$ be the distribution function of $X$, $l(x;\theta) = \log f(x;\theta)$ and write $F_s(x;\theta) = \partial_s F(x;\theta)$, $l_s(x;\theta) = \partial_s \log f(x;\theta)$. Define the functions

$$(3.2) \qquad h_r = g^{st} \int (F_s l_r - F_r l_s)\frac{\partial l_t}{\partial x}dx,$$

where the integration is over the (common) support of $F(x;\theta)$. Finally write $\lambda^J = \log \pi^J$, where $\pi^J(\theta) \propto |g(\theta)|^{1/2}$ is Jeffreys' prior.

**Theorem 3.2.** *Suppose that there exists a UPMP, $\pi$, for quantile matching. Then $\pi$ is the unique UPMP and the partial derivatives of $\lambda = \log \pi$ are given by*

$$(3.3) \qquad \partial_r \lambda(\theta) = \partial_r \lambda^J(\theta) + h_r(\theta).$$

*Proof.* We begin by expressing $g(\theta)$ in terms of the functions $\mu_t(\theta;\alpha)$ defined at (2.3). By differentiation of equation (2.2) with respect to $\alpha$ we see that $-f(q;\theta)\partial q/\partial \alpha = 1$, while differentiation of equation (2.3) gives

$$(3.4) \qquad \partial\mu_j(\theta,\alpha)/\partial\alpha = -f_j(q;\theta)\partial q/\partial\alpha = l_j(q;\theta),$$

on substitution of $\partial q/\partial \alpha$ from the previous relation. It follows that

$$(3.5) \quad g_{ij}(\theta) = \int l_i(q;\theta) l_j(q;\theta) f(q;\theta) dq = \int_0^1 \left(\frac{\partial \mu_i(\theta,\alpha)}{\partial \alpha}\right)\left(\frac{\partial \mu_j(\theta,\alpha)}{\partial \alpha}\right) d\alpha.$$

Suppose that there exists a UPMP $\pi$. Differentiation of equation (3.1) with respect to $\alpha$ and multiplication by $\partial \mu_r/\partial \alpha$ gives the equation

$$g^{st}\frac{\partial \mu_t}{\partial \alpha}\frac{\partial \mu_r}{\partial \alpha}\partial_s \lambda + \frac{\partial \mu_r}{\partial \alpha}\partial_s\left\{g^{st}\frac{\partial \mu_t}{\partial \alpha}\right\} = 0.$$

Since this relation must hold for all $0 < \alpha < 1$, integration over $0 < \alpha < 1$ gives

$$(3.6) \quad g^{st}\left\{\int_0^1 \frac{\partial \mu_t}{\partial \alpha}\frac{\partial \mu_r}{\partial \alpha} d\alpha\right\}\partial_s\lambda + \int_0^1 \frac{\partial \mu_r}{\partial \alpha}\partial_s(g^{st}\frac{\partial \mu_t}{\partial \alpha})d\alpha = 0.$$

But from (3.5) the left-hand side of (3.6) is $g^{st}g_{tr}\partial_s\lambda = \delta_r^s \partial_s\lambda = \partial_r\lambda$, where $\delta_r^s$ is the Kronecker delta function. Also, since $\partial_s(g^{st}g_{tr}) = \partial_s(\delta_r^s) = 0$, the product rule gives

$$0 = g^{st}\int_0^1 \partial_s\left(\frac{\partial \mu_r}{\partial \alpha}\right)\frac{\partial \mu_t}{\partial \alpha}d\alpha + \int_0^1 \frac{\partial \mu_r}{\partial \alpha}\partial_s\left(g^{st}\frac{\partial \mu_t}{\partial \alpha}\right)d\alpha$$

so that (3.6) becomes

$$\partial_r\lambda = g^{st}\int_0^1 \partial_s\left(\frac{\partial \mu_r}{\partial \alpha}\right)\frac{\partial \mu_t}{\partial \alpha}d\alpha.$$

This expression gives the partial derivatives of $\lambda = \log \pi$ and, furthermore, establishes that $\pi$ is the unique UPMP. We now show that this expression is equivalent to (3.3).

We first obtain the partial derivatives of $\lambda^J$. From a standard result for the derivative of a matrix determinant, we have

$$\begin{aligned}\partial_r \lambda^J &= \frac{1}{2}\partial_r \log|g| = \frac{1}{2}g^{st}\partial_r g_{st}\\
&= \frac{1}{2}g^{st}\partial_r \int_0^1 \left(\frac{\partial \mu_s}{\partial \alpha}\right)\left(\frac{\partial \mu_t}{\partial \alpha}\right) d\alpha = g^{st}\int_0^1 \partial_r\left(\frac{\partial \mu_s}{\partial \alpha}\right)\left(\frac{\partial \mu_t}{\partial \alpha}\right) d\alpha,\end{aligned}$$

again using (3.5). The difference between the $r$th partial derivatives of $\lambda$ and $\lambda^J$ is therefore

$$(3.7) \quad \partial_r\lambda - \partial_r\lambda^J = g^{st}\int_0^1 \frac{\partial}{\partial \alpha}(\partial_r\mu_s - \partial_s\mu_r)\frac{\partial \mu_t}{\partial \alpha}d\alpha.$$

Differentiation of (2.2) with respect to $\theta_r$ gives, writing $q = q(\theta,\alpha)$, $q_r f(q;\theta) + \mu_r(\theta,\alpha) = 0$, from which we obtain

$$\partial_r \mu_s(\theta,\alpha) = \int_q^\infty f_{rs}(u;\theta)du - f_s(q;\theta)q_r = \int_q^\infty f_{rs}(u;\theta)du - l_s(q;\theta)\mu_r(\theta,\alpha).$$

Furthermore, we have

$$\frac{\partial}{\partial \alpha}\{l_s(q;\theta)\mu_r(\theta,\alpha)\} = \frac{\partial q}{\partial \alpha}\frac{\partial l_s(q;\theta)}{\partial q}\mu_r(\theta,\alpha) + l_s(q;\theta)l_r(q;\theta).$$

It now follows from these two relations that

$$\frac{\partial}{\partial \alpha}(\partial_r \mu_s - \partial_s \mu_r) = \frac{\partial q}{\partial \alpha}\left(\frac{\partial l_r(q;\theta)}{\partial q}\mu_s - \frac{\partial l_s(q;\theta)}{\partial q}\mu_r\right).$$

Substituting into equation (3.7) gives

$$(3.8) \qquad h_r = g^{st}\int\left(F_s\frac{\partial l_r}{\partial q} - F_r\frac{\partial l_s}{\partial q}\right)l_t dq$$

on the change of variables from $\alpha$ to $q$, using equation (3.4) and on noting that $\mu_s(\theta, \alpha(q,\theta)) = -F_s(q;\theta)$. Next note that the indefinite integral

$$\int F_s(q;\theta)\frac{\partial l_r(q;\theta)}{\partial q}dq = F_s(q;\theta)l_r(q;\theta) - \int l_s(q;\theta)l_r(q;\theta)dq$$

from which it follows by an integration by parts that (3.8) is equivalent to (3.2), as required. □

In the case $p=1$ we have $h_r = 0$ so the unique UPMP is Jeffreys' prior, as given in Theorem 1 of [4]. For the location-scale model $f(x;\theta) = \theta_2^{-1}f^*(\theta_2^{-1}(x-\theta_1))$ discussed in Section 1, it can be verified that the solution to (3.3) is $\pi(\theta) \propto \theta_2^{-1}$, which is the right Haar prior for this model under the location-scale transformation group. In general a necessary condition for there to be a UPMP is that $h_r$ be a derivative field. The condition is not sufficient, however, as Jeffreys' prior always satisfies equation (3.3) in the case $p=1$ but we know from [4] that Jeffreys' prior is not necessarily a UPMP. When $p > 1$ the condition that $h_r$ be a derivative field is a very strong one when the model is not transformational. We have been unable to construct a two-dimensional example that is not transformational and that satisfies this condition. Even given a model satisfying this condition, the resulting prior may still not satisfy (2.5) for all $\alpha$. Thus it would seem that UPMPs rarely exist outside of transformation models. The major point of Theorem 3.2, however, is to show that if a UPMP does exist then it is unique.

Note that, whether or not a UPMP exists, when $h_r$ is a derivative field then (3.2) defines a unique prior $\pi$ which, from the proof of the Theorem 3.2, satisfies the equation $\int_0^1 \left(\frac{\partial \mu_r}{\partial \alpha}\right)\left(\frac{\partial \epsilon}{\partial \alpha}\right) d\alpha = 0$, where $\epsilon(\theta, \alpha)$ is the $O(n^{-1})$ error term (2.4). Assuming that $\partial \mu_r/\partial \alpha$ is well behaved at $\alpha = 0$ and $\alpha = 1$, integration by parts shows that this is equivalent to

$$\int_0^1 \frac{\partial^2 \mu_r}{\partial \alpha^2}\epsilon \, d\alpha = 0$$

for all $\theta$ and $r$. These relations give some sort of average prediction error, but it is unclear what precise interpretation can be given to them.

Finally, when there exists a suitable group structure as discussed earlier then we know that $\partial_i \lambda$ must be $\partial_i \lambda^H$. Furthermore, since Jeffreys' prior is the left Haar prior, it follows that $h_i(\theta) = \partial_i \log \Delta(\theta^{-1})$, where $\Delta$ is the modulus of $\Omega$ and $\theta^{-1}$ is the group inverse of $\theta$.

## 4. UPMPs: highest predictive density region matching

We consider now the case where $X$ is possibly vector-valued. The question of the existence of UPMPs for highest predictive density regions in this case is not so straightforward as the quantile matching case discussed in Section 3. In particular, if there exists a suitable group structure, as in Section 3, since the prediction region $H(\pi, \alpha, d)$ defined in Section 1 is not invariant under transformation of $X$ (unless the group is affine), the associated right Haar prior is not necessarily a UPMP. We also know that when a UPMP does exist it may not be unique. This was illustrated in Example 4 in [4] of the bivariate normal model with unknown covariance matrix; we will return to this example in Example 1 below.

The scalar parameter case is straightforward, however. For each $\alpha$ the prior

$$\pi(\theta) \propto g(\theta)\{\xi_1(\theta, \alpha)\}^{-1}$$

is the unique solution to (2.6). It follows that there exists a UPMP prior if and only if $\xi_1(\theta, \alpha) = Q(\theta)R(\alpha)$, as was noted in [4] where examples are given in which this condition does and does not hold. Unlike the case of quantile matching, however, the unique solution when it exists is not necessarily Jeffreys' prior. For example, in [4] it is shown that a unique UPMP exists for the $N(\theta, \theta)$ model but this is not Jeffreys' prior.

The multiparameter case is more difficult. The simplest situation is when $\xi_t(\theta, \alpha)$ is of the form

(4.1) $$\xi_t(\theta, \alpha) = Q_t(\theta)R(\alpha).$$

Then every UPMP will be a solution of the Lagrange PDE

(4.2) $$\partial_s\{g^{st}(\theta)Q_t(\theta)\pi(\theta)\} = 0.$$

This equation may have no solutions or an infinite number of solutions.

**Example 1.** Consider the bivariate normal model with zero means and unknown standard deviations $\sigma_1, \sigma_2$ and correlation coefficient $\rho$. Let $\Sigma$ be the covariance matrix of $X$. We work with the orthogonal parameterisation

$$T^{-1} = \begin{pmatrix} \theta_1 & 0 \\ \theta_2\theta_3 & \theta_2 \end{pmatrix},$$

where $\Sigma = TT'$ and $T$ is the left Cholesky square root of $\Sigma$. It can then be shown that the information matrix is

$$g(\theta) = \mathrm{diag}(2\theta_1^{-2}, 2\theta_2^{-2}, \theta_1^{-2}\theta_2^2).$$

Furthermore, by transforming to $Z = T^{-1}X$, it can be shown that

$$m(\theta, \alpha) = \theta_1\theta_2(1-\alpha)/(2\pi), \; \xi_1(\theta, \alpha) = \theta_1^{-1}R(\alpha), \; \xi_2(\theta, \alpha) = \theta_2^{-1}R(\alpha), \; \xi_3(\theta, \alpha) = 0,$$

where $R(\alpha) = -(1-\alpha)\log(1-\alpha)$. Thus $\xi_t(\theta, \alpha)$ is of the form (4.1). Therefore the UPMP priors are all the solutions of the PDE (4.2) with $Q_1(\theta) = \theta_1^{-1}$, $Q_2(\theta) = \theta_2^{-1}$ and $Q_3(\theta) = 0$. The general solution is found to be

$$\pi(\theta) \propto \theta_1^{-2}h(\theta_2^{-1}\theta_1, \theta_3),$$

where $h$ is an arbitrary positive function. Notice that the leading term $\theta_1^{-2}$ is $|g(\theta)|^{1/2}$, so Jeffreys' prior is a UPMP. In terms of $(\sigma_1, \sigma_2, \rho)$ we have

$$\theta_1 = \sigma_1^{-1}, \; \theta_2 = \sigma_2^{-1}(1-\rho^2)^{-1/2}, \; \theta_3 = -\rho\sigma_1^{-1}\sigma_2$$

with Jacobian of transformation $\sigma_1^{-3}\sigma_2^{-1}(1-\rho^2)^{-3/2}$. With suitable re-expression of $h$ we find that

$$\pi(\sigma_1, \sigma_2, \rho) \propto \pi^J(\sigma_1, \sigma_2, \rho) H(\sigma_1^{-1}\sigma_2, (1-\rho^2)^{1/2}), \tag{4.3}$$

where $\pi^J(\sigma_1, \sigma_2, \rho) \propto \sigma_1^{-1}\sigma_2^{-1}(1-\rho^2)^{-3/2}$ is Jeffreys' prior and $H$ is an arbitrary positive function. This is a very wide class of priors. In particular, taking $h(x, y) = x^a y^b$, we see that all priors of the form $\pi^J(\sigma_1, \sigma_2, \rho)(\sigma_1^{-1}\sigma_2)^a(1-\rho^2)^b$ are UPMPs. Taking $a = 1, b = 1/2$ we obtain $\sigma_1^{-2}(1-\rho^2)^{-1}$, which can be shown to be the right Haar prior arising from the group of transformations $T^{-1}X$ on the sample space, where $T$ is a lower triangular matrix with positive diagonal elements. This group is isomorphic to $\Omega$ and since in this case the region $A$ is invariant it follows from [13] that this prior must be a UPMP. Similarly, all right Haar priors arising from transformations of the form $T^{-1}MX$, with $M$ a fixed non-singular matrix, are included in (4.3).

We now return to the general analysis of equation (2.6). In Theorem 4.1 below, when we say that the functions $\xi_t(\theta, \alpha)$ are *linearly independent* we shall mean that they are linearly independent as functions of $\alpha$ for fixed $\theta$.

**Theorem 4.1.** *Suppose that the functions $\xi_t(\theta, \alpha)$ are linearly independent and that there exists a UPMP, $\pi$, for highest predictive density region matching. Then $\pi$ is the unique UPMP and the partial derivatives of $\lambda = \log \pi$ are given by*

$$\partial_j \lambda = -b^{ri}g_{ij} \int_0^1 \frac{\partial \xi_r}{\partial \alpha} \partial_s \left( g^{st} \frac{\partial \xi_t}{\partial \alpha} \right) d\alpha, \tag{4.4}$$

*where $(b^{ri}(\theta))$ is the inverse of the non-singular matrix function $(b_{ij}(\theta))$ with $(i,j)$th element*

$$b_{ij}(\theta) = \int_0^1 \left( \frac{\partial \xi_i(\theta, \alpha)}{\partial \alpha} \right) \left( \frac{\partial \xi_j(\theta, \alpha)}{\partial \alpha} \right) d\alpha. \tag{4.5}$$

*Proof.* We begin by showing that the matrix $(b_{ij}(\theta))$ is non-singular for all $\theta \in \Omega$ if and only if the functions $\xi_t(\theta, \alpha)$ are linearly independent. From the definition (4.5), we see that in general $(b_{ij}(\theta))$ is positive semidefinite and is therefore singular for all $\theta \in \Omega$ if and only if, for each $\theta$, there exist functions $x^t(\theta)$, not all zero, for which $b_{ij}(\theta)x^i(\theta)x^j(\theta) = 0$. This is equivalent to the condition

$$\int_0^1 \left( \frac{\partial x^t(\theta) \xi_t(\theta, \alpha)}{\partial \alpha} \right)^2 d\alpha = 0,$$

which in turn holds if and only if $\partial(x^t(\theta)\xi_t(\theta, \alpha))/\partial \alpha = 0$ for all $\theta$ and $\alpha$. Since $\xi_t(\theta, 1) = 0$ it follows that a necessary and sufficient condition for the singularity of $(b_{ij}(\theta))$ is the existence of $x^t(\theta)$, not all zero, such that $x^t(\theta)\xi_t(\theta, \alpha) = 0$ for all $\theta$ and $\alpha$. That is, the functions $\xi_t(\theta, \alpha)$ are linearly dependent.

We now apply the product rule to (2.6) to give equation (3.1) with $\mu_t$ replaced by $\xi_t$. Exactly as in the proof of Theorem 3.2, we differentiate this equation with respect to $\alpha$, multiply by $\partial \xi_r/\partial \alpha$ and integrate over $0 < \alpha < 1$ to give

$$g^{st}b_{tr}\partial_s \lambda + \int_0^1 \frac{\partial \xi_r}{\partial \alpha} \partial_s \left( g^{st} \frac{\partial \xi_t}{\partial \alpha} \right) d\alpha = 0.$$

Finally, under the condition of the theorem the matrix $(b_{ij}(\theta))$ is non-singular and equation (4.4) follows on multiplying both sides of the above expression by $b^{ri}g_{ij}$. □

In the case $p = 1$ we know that a UPMP exists if and only if $\xi_1(\theta, \alpha) = Q(\theta)R(\alpha)$, in which case

$$b_{11}(\theta) = \{Q(\theta)\}^2 \int_0^1 \{R(\alpha)\}^2 d\alpha.$$

Equation (4.4) then becomes $d\lambda/d\theta = d\log(g^{-1}\theta)/d\theta$, giving $\pi(\theta) \propto \{Q(\theta)\}^{-1}g(\theta)$ in agreement with the earlier discussion. In the multiparameter case, unlike Theorem 3.2, there does not appear to be any simple further development of (4.4). Returning to the univariate location-scale model $f(x; \theta) = \theta_2^{-1} f^*(\theta_2^{-1}(x - \theta_1))$, it can be verified that the functions $\xi_t(\theta, \alpha)$ are linearly independent and, as in Section 3, that the right Haar prior $\pi(\theta) \propto \theta_2^{-1}$ under the location-scale transformation group is the solution to (4.4). When $p > 1$ the condition that the right-hand side of (4.4) be a derivative field is very strong when the model is not transformational and we have been unable to find a two-dimensional example that is not a transformation model satisfying this condition. Again, as in Section 3, even for such an example the resulting prior may still not satisfy (2.6) for all $\alpha$. Thus it would seem that unique UPMPs rarely exist outside of transformation models. As with Theorem 3.2, the major point of Theorem 4.1 is to show that, under the conditions of the Theorem, if a UPMP does exist then it is unique.

Note that when $p > 1$ Theorem 4.1 does not apply to the case (4.1) since the functions $\xi_t(\theta, \alpha)$ are linearly dependent and hence the matrix $(b_{ij}(\theta))$ is singular. A more general sufficient condition for linear dependence of the $\xi_t(\theta, \alpha)$ is

(4.6) $$\xi_t(\theta, \alpha) = U_t(\theta)S(\theta, \alpha).$$

Note that this is also a necessary condition for linear dependence in the case $p = 2$.

Suppose that (4.6) holds and that there exists a UPMP $\pi$. Then from equation (2.6) we see that

$$g^{st}U_t \partial_s \lambda + g^{st}U_t \partial_s \log S + \partial_s(g^{st}U_t) = 0$$

for all $\alpha$, which implies that the function $g^{st}(\theta)U_t(\theta)\partial_s \log S(\theta, \alpha)$ must be free from $\alpha$. Since no boundary conditions are imposed on the solutions to the resulting Lagrangian PDE, it follows that $\pi$ must be one of an infinite number of solutions. Thus, under condition (4.6), either there is no UPMP or there is an infinite number of UPMPs. Note that (4.1) is a special case of (4.6).

It might appear at first sight that it is also possible to have an infinite number of UPMPs in the case of quantile matching, which would contradict the result of Theorem 1. However, using a parallel argument to that given above, we see that the structure (4.6) for $\mu_t(\theta, \alpha)$ cannot occur, as this would imply singularity of Fisher's information matrix.

Finally, the case

(4.7) $$\xi_t(\theta, \alpha) = Q_t(\theta)R_t(\alpha),$$

which is a generalisation of the simple case (4.1), is of some interest. It is easily seen that in this case the linear independence of the functions $\xi_t(\theta, \alpha)$ is equivalent to the linear independence of the functions $R_t(\alpha)$. Furthermore, the matrix $(b_{ij}(\theta))$ will be non-singular for all $\theta$ if and only if the matrix with $(i, j)$th element

$$a_{ij} = \int_0^1 \left(\frac{\partial R_i}{\partial \alpha}\right)\left(\frac{\partial R_j}{\partial \alpha}\right) d\alpha$$

is positive definite. This turns out to be the case for the location-scale models discussed earlier.

**Example 2.** Consider the multivariate location model with density $f(x;\theta) = f^*(x_1 - \theta_1, \ldots, x_p - \theta_p)$. The region $A$ here is invariant under the group of transformations $x + a$, $a \in \mathcal{R}^p$ and it follows from [13] that the right Haar prior is an exact UPMP. Here the right Haar prior is also Jeffreys' prior, both being constant. We now investigate conditions under which this is the unique UPMP. As in [4], we find that $m(\theta, \alpha) = m(\alpha)$, free from $\theta$, and $\xi_t(\theta, \alpha) = R_t(\alpha)$, which is of the form (4.7) with $Q_t(\theta) = 1$ for all $t$. It follows from the above discussion that Jeffreys' prior is the unique UPMP if and only if the functions $R_t(\alpha)$ are linearly independent. In that case, since both $g^{st}$ and $\xi_t$ are free from $\theta$, the right-hand side of (4.4) is zero and, again, the unique UPMP is the uniform prior.

For many standard models, however, the functions $R_t(\alpha)$ will be linearly dependent. Suppose, for example, that $f^*$ is elliptically symmetric, so that $f^*(z) = H(z'Cz)$ for some positive definite matrix $C$. Then it can be checked that $R_t(\alpha) = Q_t R(\alpha)$, which is of the form (4.1) with $Q_t$ free from $\theta$. The functions $R_t(\alpha)$ are clearly linearly dependent and hence, since we know that there exists at least one UPMP, there will be an infinite number of UPMPs. For example, in the case where $f^*$ is spherically symmetric, we have $Q_t = Q$ and the Lagrange PDE (4.2) becomes $\sum_s \partial_s \lambda = 0$. The solutions of this equation are of the form $\pi(\theta) \propto \exp\{h(\theta_2 - \theta_1, \ldots \theta_p - \theta_1)\}$, where $h$ is an arbitrary function. In particular, all priors of the form $\pi(\theta) \propto \exp(\sum_i a_i \theta_i)$ with $\sum_i a_i = 0$ will be uniformly matching in this case.

A similar analysis may be carried out for the multiparameter location-scale model with different location parameters, as described in [4]. Whether or not the scale parameters are assumed to be equal, there is an appropriate group of transformations for which the corresponding right Haar prior will be a UPMP. In either case $\xi_t(\theta, \alpha)$ is again of the form (4.7) so that whether or not the right Haar prior is the unique UPMP will depend on the linear independence or otherwise of the functions $R_t(\alpha)$.

When the model has no suitable group structure, we conjecture that the functions $\xi_r(\theta, \alpha)$ will always be linearly independent. To see the plausibility of this, note that the $\xi_r(\theta, \alpha)$ are linearly dependent if and only if there exist functions $x^t(\theta)$, not all zero, such that $\int_A \{x^t(\theta) l_t(x; \theta)\} f(x; \theta) dx = 0$ for all $\theta$ and $\alpha$. Since the density $f(x;\theta)$ cannot be standardised by transformation, the only way that this would seem to be possible is if $x^t(\theta) l_t(x; \theta) = 0$ for all theta. However, it is easily seen by partial differentiation w.r.t $\theta_s$ that this condition leads to $g$ being singular. This analysis therefore suggests that if the model is not transformational then there will either be no UPMP or a unique UPMP, which is then given by (4.4).

## 5. Discussion

Although it is known that exact matching of invariant prediction regions is achieved by the right Haar prior under a suitable group structure on the model, we have seen in Section 3 that there can be other priors that achieve approximate uniform predictive quantile matching, and that uniformly matching priors can exist when there is no suitable group structure, although these are rare. In common with other work on probability matching priors, predictive matching priors arise as solutions to a particular PDE, which in general can be very difficult to solve. However, in the case of uniform quantile matching, if a UPMP exists then it is unique and explicit formulae for its partial derivatives are available from Theorem 3.2.

Except in special cases, derivation of the UPMP for quantile matching via equation (3.3), or even verifying that the derivatives in (3.3) are consistent, will be intractable. An attractive alternative would be to use a data-dependent approximation of the UPMP based on a local prior of the form

$$\partial_r \lambda(\theta, \theta_0) = \partial_r \lambda^J(\theta) + h_r(\theta_0).$$

See [14] for a derivation of data-dependent matching priors for marginal posterior distributions. Furthermore, since a data-dependent prior of this form will always exist, there may be cases for which it will be uniformly matching even when there is no $\alpha$-free solution of (2.5). Although the posterior distribution arising from such a prior would not always have a strict Bayesian interpretation, use of the corresponding predictive distribution could provide a useful mechanism for constructing frequentist prediction regions with good coverage properties. It would be of interest to conduct simulation experiments in order to assess the predictive coverage afforded by such priors.

The case of highest predictive density regions is more complex. As discussed in Section 4, there will either be a unique solution or else there will be infinitely many solutions, depending on the linear independence or otherwise of the functions $\xi_r(\theta, \alpha)$. Thus in any particular example it is necessary to examine carefully the structure of the functions $\xi_r(\theta, \alpha)$. If the statistical model has a suitable group structure then this task is usually eased. One could also investigate local priors when the matrix $(b_{ij}(\theta))$ is invertible.

In the case of univariate observations, the results provide some guidance on the choice of objective prior if the main goal is to carry out Bayesian prediction and low predictive coverage probability bias is desired. In relation to the determination of an objective prior, for multivariate data the situation is less clear. When the functions $\xi_r(\theta, \alpha)$ are linearly dependent, as often occurs in transformation models, there will usually be an infinite number of UPMPs. Thus other considerations will need to be invoked in order to narrow down the choice of prior. For example, one might consider priors that are simultaneously predictive and posterior probability matching, reference priors ([2]) or priors that are minimax under suitable decision rules; in particular, for minimax prediction loss see, for example, [12] and [15].

**Appendix: Proof of Theorem 3.1**

*Proof.* Let $a \in \Omega$ and consider the transformation $\phi = a\theta$. Let $J(\theta, a) = \partial \phi / \partial a$ be the Jacobian matrix of this transformation for fixed $\theta$. Then the right Haar prior is $\pi^H(\theta) \propto |J(\theta, e)|^{-1}$, where $|J(\theta, a)|$ is the determinant of $J(\theta, a)$; see, for example, [1].

Write $\tilde{\phi}_s^r(\theta, a) = \partial \phi_r(\theta, a)/\partial a_s$ and define $\tilde{\phi}_s^r(\theta) = \tilde{\phi}_s^r(\theta, e)$, where $e$ is the identity element of the group, so that $\pi^H(\theta) \propto |(\tilde{\phi}_s^r(\theta))|^{-1}$. Finally, let $(a_r^s(\theta))$ be the matrix inverse of $(\tilde{\phi}_s^r(\theta))$. A standard result for the derivative of a matrix determinant then gives

(5.1) $$\partial_s \lambda^H(\theta) = -a_r^u(\theta) \partial_s \tilde{\phi}_u^r(\theta),$$

where $\lambda^H(\theta) = \log \pi^H(\theta)$.

Define $\phi_s^r(\theta, a) = \partial_s \phi_r = \partial \phi_r / \partial \theta_s$, with matrix inverse $\theta_r^s(\theta, a) = \partial \theta_s / \partial \phi_r$. Since the definition of the right Haar prior depends on a specific group of transformations on the parameter space, it is natural to regard Fisher's information as a Riemannian metric tensor associated with the differentiable manifold of probability densities $f(\cdot; \theta)$, $\theta \in \Omega$. This facilitates the study of the transformational properties of the quanitities $g^{st}(\theta)$ and $\mu_t(\theta,\alpha)$ appearing in the PDE (3.1). First, from the invariance of the problem under $G$ and the contravariant tensorial property of $g^{st}$, we have $\bar{g}^{ij}(\phi) = g^{st}(\theta)\phi_s^i\phi_t^j$, where $\bar{g}^{ij}$ is the inverse Fisher information in the $\phi$-parameterisation. Again using the invariance properties, it is seen that $\bar{\mu}_j(\phi,\alpha) = \mu_k(\theta,\alpha)\theta_j^k$, where $\bar{\mu}_j(\phi,\alpha)$ is the function (2.3) in the $\phi$-parameterisation. Now write $u^s(\theta,\alpha) = g^{st}(\theta)\mu_t(\theta,\alpha)$ and $\bar{u}^s(\phi,\alpha) = \bar{g}^{st}(\phi)\bar{\mu}_t(\phi,\alpha)$. Then

$$
\begin{aligned}
\bar{u}^i(\phi,\alpha) &= g^{st}(\theta)\mu_k(\theta,\alpha)\phi_s^i\phi_t^j\theta_j^k \\
&= g^{st}(\theta)\mu_k(\theta,\alpha)\phi_s^i\delta_t^k = u^s(\theta,\alpha)\phi_s^i,
\end{aligned}
\tag{5.2}
$$

where $\delta_t^k$ is the Kronecker delta function.

Now differentiate both sides of (5.2) with respect to $a_r$ to give

$$
\partial_s\bar{u}^i(\phi,\alpha)\tilde{\phi}_s^r(\theta,a) = u^s(\theta,\alpha)\partial\phi_s^i(\theta,a)/\partial a_r = u^s(\theta,a)\partial_s\phi_r^i(\theta,a).
\tag{5.3}
$$

Finally, setting $a = e$ and multiplying both sides of (5.3) by $a_i^r(\theta)$ gives

$$
\partial_s u^i(\theta,\alpha)a_i^r(\theta)\tilde{\phi}_r^s(\theta) = u^s(\theta,\alpha)a_i^r(\theta)\partial_s\tilde{\phi}_r^i(\theta).
\tag{5.4}
$$

Since $(a_r^s(\theta))$ is the matrix inverse of $(\tilde{\phi}_s^r(\theta))$, the left-hand side of (5.4) is $\partial_s u^i(\theta,\alpha)\delta_i^s = \partial_s u^s(\theta,\alpha)$, whereas the right-hand side is $-u^s(\theta,\alpha)\partial_s\lambda^H(\theta)$ from (5.1). It follows that the right Haar prior $\pi^H$ is a solution of equation (3.1) and hence of equation (2.5). □